\documentclass[a4paper,11pt]{article}

\newtheorem{theorem}{Theorem}
\newtheorem{proposition}[theorem]{Proposition}
\newtheorem{lemma}[theorem]{Lemma}
\newtheorem{corollary}[theorem]{Corollary}
\newtheorem{definition}[theorem]{Definition}

\newtheorem{example}[theorem]{Example}

\newtheorem{remark}[theorem]{Remark}

\newcommand{\qed}{\nobreak \ifvmode \relax \else
      \ifdim\lastskip<1.5em \hskip-\lastskip
      \hskip1.5em plus0em minus0.5em \fi \nobreak
      \vrule height0.75em width0.5em depth0.25em\fi}
\newcommand{\dimn}{n}
\newcommand{\Tor}{{\mathrm{Tor}}}

\makeindex % please use the style svind.ist with your makeindex program
\usepackage{makeidx}     % allows index generation
\usepackage{cite}        % for citations
\usepackage{amsfonts}
\usepackage{amssymb}
\usepackage{amscd}
\usepackage{epsfig}
\usepackage{amsrefs}
\usepackage{url}
\def\C{\mathbf{C}}
\def\F{\mathbf{F}}

\def\R{\mathbf{R}}

\def\rmi{\uppercase\expandafter{\romannumeral1}}
\def\rmii{\uppercase\expandafter{\romannumeral2}}
\def\rmiii{\uppercase\expandafter{\romannumeral3}}
\def\rmv{\uppercase\expandafter{\romannumeral5}}
\def\rmvi{\uppercase\expandafter{\romannumeral6}}
\def\rmvii{\uppercase\expandafter{\romannumeral7}}
\def\rmviii{\uppercase\expandafter{\romannumeral8}}

\def\FF{\mathcal{F}}          % functions on something

\def\X{\mathcal{X}}              % ham vector field
\def\can#1{\left\langle#1\right\rangle}

\def\pb#1{\left\{#1\right\}}

\def\lb#1{\[#1\]}
\def\LB{[\cdot\,,\cdot]}

\def\({\left(}
\def\){\right)}
\def\[{\left[}
\def\]{\right]}

\def\diff{\mathsf{d}}

\newenvironment{eqn*}[1][1.5]
  {$$\renewcommand{\arraystretch}{#1}
      \begin{array}{rcl}}
      {\end{array}$$}

\def\comment#1{}  % to comment out temporarily
\def\p{\partial}
\def\pp#1#2{\frac{\p #1}{\p #2}}
\begin{document} 
\nocite{*}

\title{A Lagrangian for Hamiltonian vector fields on singular Poisson manifolds}
\author{Yahya Turki
\thanks{Universit\'e de Lorraine, 
Institut Elie Cartan de Lorraine
UMR 7502, Metz, F-57045, France.
\texttt{yahya.turki@univ-lorraine.fr} 
}  
\thanks{Universit\'e de Monastir, Facult\'e des Science de Monastir, Avenue de l'environnement 5019 Monastir, Tunisie  
}
}

\maketitle
\begin{abstract}
On a manifold equipped with a bivector field, we introduce for every Hamiltonian a Lagrangian 
on paths valued in the cotangent space whose stationary points projects onto Hamiltonian vector fields.
We show that the remaining components of those stationary points tell whether the bivector field is Poisson
or at least defines an integrable distribution - a class of bivector fields generalizing twisted Poisson structures 
that we study in detail.
\end{abstract}

\tableofcontents
 % the title, authors, abstract and all that stuff
\newpage
\section{Introduction}

Poisson $\sigma$-models, developed by Ikeda \cite{Ikeda}
and Shaller-Strobl \cite{SchallerStrobl} is now a well-developed theory, 
well-known to give an alternative approach to Kontsevich star-product, see Cattaneo-Felder
\cite{CattaneoFelder}. Poisson $\sigma$-models is based on the study of a certain functional, 
defined on vector bundle morphisms from the tangent space of a cylinder to the cotangent space of a manifold $M$
equipped with a bivector field $\pi$. Our purpose is to study a very natural functional ${\mathcal L}^H$,
which is in the same spirit as the one defining Poisson $\sigma$-models, but 
the dimension of the source manifold is $1$ and $M$ comes equipped with a bivector field
and a function. Explicitly, this functional is given by:
\begin{equation}
   \label{eq:defL0}
 {\mathcal L}^H (\alpha) = \int_0^1 \left\langle \left. {\mathcal X}_H \right|_{x(t)}  -  \frac{\diff x(t)}{\diff t} , \alpha(t) \right\rangle \diff t ,
  \end{equation}
  with $\X_H$ being the Hamiltonian vector field.
Our functional is also inspired by the celebrated functional of Weinstein \cite{Weinstein}, which is defined 
on paths valued in an exact symplectic manifold $(M,\omega=\diff \beta)$:
\begin{equation}
   \label{eq:defW}
 {\mathcal L}^H (\alpha) = \int_0^1 \left( \left\langle \beta_{x(t)} , \frac{\diff x(t)}{\diff t}  \right\rangle - H(x(t)) \right) \diff t,
 \end{equation}
 See section \ref{sec:relation} for a more precise relation between all those functional.
The initial purpose of the present article is to state results of the form "A bivector field is of type X if and only if 
the stationary points of ${\mathcal L}^H$ are, for all $H$, of type Y". We were especially interested to find the Y corresponding to X = "Poisson" 
and the X corresponding to Y = "cotangent paths". To be able to state such a result, we are obliged to introduce several 
new notions.
\begin{enumerate}
 \item Quasi-cotangent paths (with respect to a function $H \in C^\infty (M)$), which are simply paths $\alpha(t)$ in $T^*M$ for which $\pi^\#(\alpha)-\frac{\diff x(t)}{\diff t}$ 
 is $\X_H$-invariant, with $x$ the base path of $M$. In particular cotangent paths are quasi-cotangent.
 \item Foliated bivector fields, which are defined to be those for which the distribution $\pi^\#(T^*M)$ is integrable. This is a subtle notion that
 we claim to have an interest of its own. For every singular integrable distribution in the sense of Sussmann \cites{Stefan,Sussmann},
 it is natural to look for an algebroid whose image through the anchor map gives the distribution, at least locally \cite{androulidakis}.
 For foliated bivector fields, this algebroid is not a priori given and 
 does not seem easy to guess. This open question will be addressed in a subsequent paper.
  \item Weakly foliated bivector fields, which are defined to be those for which the distribution $\pi^\#(T^*M)$ is integrable, so to say,
  at every point, i.e. such that $\left. \left[\pi^\#(T^*M),\pi^\#(T^*M)\right] \right|_m$ is the image of $\pi^\#_m$ at every point $m$ of the manifold $M$. 
\end{enumerate}
Our main result states the following, see theorem \ref{theo:caract_chemins_cotangents} for more details.
Let $M$ be a manifold and $\pi $ be a bivector field on $M$.
\begin{enumerate}
\item  If X = "foliated" then Y = "cotangent paths" for every function $H \in C^\infty (M) $.
\item  If X = "weakly foliated" then  Y = "cotangent paths" for every function $H \in C^\infty (M) $. 
\item  X = "Poisson" if and only if Y = "quasi-cotangent paths" for every function $H \in C^\infty (M) $.
\end{enumerate}
The reader may notice that we do not work with loops, but with paths, which is highly surprising since from Weinstein's Lagrangian (\ref{eq:defW}),
it is periodicity that makes things interesting.
We are in fact also interested by the periodic case, but this should be postponed to a subsequent article.
%Here, we just will give an interpretation of stationary points in the periodic case.

The paper is organized as follows. Section \ref{sec:foliated} is devoted to the study of a new type of bivector fields that we call foliated bivector fields, 
which are defined to be those bivector field $\pi$ for which the distribution $\pi^\#(T^*M)$ is integrable.
As will appear in the course of section \ref{sec:foliated}, foliated bivector fields are strongly related to twisted Poisson structures
\cites{strobelklim,KosmannSchwarzbachYvette}.
Indeed, these structures, also called Poisson structures with background, are shown in proposition \ref{prop:twistedmeansfoliated} 
to be foliated. Also, every foliated bivector field comes from a twisted Poisson structure at regular points,
see proposition \ref{prop:regularpoints}. Last, each leaf of a 
foliated Poisson structure comes equipped with a twisted Poisson structure of maximal rank, see theorem \ref{theo:localstructureFoliated}.
Conformally Poisson structures are also among examples (see proposition \ref{prop:conformally}),
indeed, in contrast with Poisson structures, foliated bivector fields are a $C^\infty(M)$-module.
Also, weak Poisson structures are defined and related to their "strong" counterpart. 

In section \ref{sec:LagrangianTheory}, the Lagrangian briefly introduced in (\ref{eq:defL0}) above is studied in details,
and the main result of this paper, theorem \ref{theo:caract_chemins_cotangents}, is stated and proved. 
Section \ref{sec:relation} explains the relation with Poisson $\sigma$-models.

I express my gratitude to the University of Monastir for two "bourses d'alternance" that I received while preparing this manuscript.

\newpage

\section{Foliated bivector fields}
\label{sec:foliated}
\subsection{Definitions and notations}

We define foliated bivector fields and indicate here some definitions, notations and results which will be needed in the sequel.
From $M$ is an arbitrary manifold of dimension $\dimn$, and $A^\bullet (M) := \sum_{k=0}^\dimn A^k(M)$ is the Gerstenhaber algebra of multivector fields,
equipped with the wedge product and the Schouten-Nijenhuis bracket, that we simply denote by $\LB$.
Here, given a vector space $E$, for all $\pi \in \wedge^{2}E$ we denote by 
$\pi^\# $ the morphism of vector bundle of $E^*  \to E$ given by for all $\xi ,\eta \in E^*$:
 \begin{equation} \label{eq:dual_and_r_sharp} 
\can{\pi^\#(\xi),\eta} = \can{\pi , \xi \wedge \eta }. 
\end{equation}
When $\can{\alpha \wedge \beta , \xi \wedge \eta  } = \can{\alpha,\xi } \can{\beta ,\eta }-\can{\alpha,\eta  } \can{\beta ,\xi }$ is the natural pairing between $\wedge ^{2}E$ and $\wedge ^{2}E^{*}$.
We recall \cites{ACP,CrainicFernandes,FernandesCranicMarius} that a \emph{Poisson structure} is a bivector field $\pi\in A^{2}(M)$ satisfying $[\pi,\pi]=0$.
This property implies that the biderivation of  $\C^{\infty}(M)$ defined by 
 $$ \pb{F,G} := \langle \pi, \diff F \wedge \diff G \rangle $$
is a Lie bracket, and that the map assigning to a function $F$ its Hamiltonian vector field $ \X_F := \pi^\#(\diff F) $
is an anti-Lie algebra morphism, i.e.
 $ \lb{\X_F,\X_G} = \X_{\lb{F,G}}.$
  for all $F,G \in \C^{\infty}(M)$. This implies that:
\begin{equation}\label{eq:integrable}
 \lb{ \pi^\# (\Omega^1(M)) , \pi^\# (\Omega^1(M)) } \subset  \pi^\# (\Omega^1(M) ). 
\end{equation}
Equation (\ref{eq:integrable}) means that the locally finitely generated 
 $\C^{\infty}(M)$-module  $\pi^\# (\Omega^1(M))$ (which is a sub-module of $ A^{1}(M)$) is a \emph{integrable distribution}. 
However, a bivector $\pi$ can verify the latter property if it is not Poisson.
This is precisely the point of the following definition:
\begin{definition}\label{def:bivecteurfeuillete}
Let $M$ be a manifold. We call \emph{foliated bivector field}
a bivector field $\pi \in A^{2}(M)$ such that the image of $\pi^\# : \Omega^1(M) \to A^{1}(M)$ is a integrable distribution. 
\end{definition}
 The following lemma is immediate.
 
 \begin{lemma}\label{lem:immediate}
A bivector field $\pi \in A^{2}(M)$ is foliated if and only if for any pair $F, G$ of functions on $ M $, 
there exists a $ 1$-form $\alpha_ {F, G} $ such that $ [\X_F,\X_G]=\pi^\# (\alpha_{F,G})$.
\end{lemma}

\subsection{Examples of foliated bivector fields}

Obviously, if $\pi \in A^{2}(M)$ is Poisson, then $ \pi $ is foliated. But there are more example.

Let $M$ be a smooth manifold. A pair $(\pi, \phi )$, where $\pi$ is a bivector field and $\phi$ is a closed 3-form, is called a \emph{twisted Poisson structure} \cites{ParkJae-Suk,strobelklim,
KosmannSchwarzbachYvette} if it satisfies:
\begin{equation}
\frac{1}{2}\lb{\pi,\pi} =  \wedge ^{3} \pi^{\sharp }(\phi ) \hbox{ and } \diff \phi=0.
\end{equation}

\begin{proposition}\label{prop:twistedmeansfoliated}
For every twisted Poisson structure $(\pi,\phi )$ the bivector field $\pi$ is foliated.
\end{proposition}
\begin{proof}
The proof is a special case of a more general phenomenon: for any Lie algebroid \cite{MackenzieKirill}, the image of the anchor is an integrable distribution.
A twisted Poisson structure \cites{SeveraWeinstein,CattaneoAlbertoPing} induces a Lie algebroid structure on $T^{*}M \to M$, with anchor map $\pi^{\sharp }$, and Lie bracket on $\Omega^1(M)$ given by:
\begin{equation}
\lb{\alpha ,\beta }={\mathcal L}_{\pi^{\sharp }(\alpha )}\beta - {\mathcal L}_{\pi^{\sharp }(\beta  )}\alpha - \diff \pi(\alpha ,\beta )+
\imath_{(\pi^{\sharp}(\alpha ) \wedge \pi^{\sharp}(\beta ) )} \phi. 
\end{equation}
In particular, for all $F,G \in \C^{\infty}(M)$, we have:
\begin{equation}\label{eq:bracketOfExactForms}
\lb{ \diff F, \diff G}= \diff \pb{F,G}+ \imath_{(\X_F ) \wedge (\X_G) }\phi,
\end{equation}
with, as before, $\pb{F,G}=\pi(\diff F,\diff G)$. Since $\pi^{\sharp }$ is the anchor map of this Lie algebroid structure, 
$ \pi^{\sharp}(\lb{ \diff F, \diff G}) =\lb{ \pi^\sharp (\diff F), \pi^\sharp (\diff G)} = \lb{\X_F,\X_G}$.
Applying $\pi^\sharp$ to both sides of (\ref{eq:bracketOfExactForms}) gives
\begin{equation}
\lb{\X_F,\X_G}=\pi^{\sharp }( \diff \pb{F,G}+  \imath_{(\X_F ) \wedge (\X_G ) }\phi )
\end{equation}
and consequently:
$$\lb{\X_F,\X_G}= \pi^\# (\alpha_{F,G}) \hbox{ where } \alpha_{F,G}=\diff \pb{F,G}+ \imath_{(\X_F ) \wedge (\X_G) }\phi .$$
The result follows from lemma \ref{lem:immediate}.
\end{proof}
Twisted Poisson structures are in fact the generic example
in the sense that at regular points, foliated bivector fields always arise 
(locally) from a twisted Poisson structure.
Recall that a bivector field $\pi^\#$ is \emph{regular} 
when the rang of $\pi_{m'}^\#: T_{m'}^{*}M \mapsto T_{m'}M $ does not depend 
on $m' \in {\mathcal U}$, where ${\mathcal U}$ is a neighborhood of $m'$, 
and \emph{regular} at a point $m' \in M$ when its restriction to a neighborhood is
regular.

\begin{proposition}\label{prop:regularpoints}
Let $\pi$ be a foliated bivector field on a manifold $M$.
\begin{enumerate}
 \item For every regular point $m\in M$ of $\pi$, 
there is a neighborhood ${\mathcal U}$ and a closed $3$-form $\varphi \in \Omega^2({\mathcal U})$,
such that $(\pi,\varphi)$ is a twisted Poisson structure on ${\mathcal U}$.
 \item Moreover, if $M$ is oriented, there exists a closed $3$-form $\varphi$ defined on the open subset
 of all regular points $ {\mathcal U}_{reg}$,
such that $(\pi,\varphi)$ is a twisted Poisson structure on ${\mathcal U}_{reg}$.
\end{enumerate}

\end{proposition}
\begin{proof}
It suffices to prove item 2, since every point has a neighborhood which is an oriented submanifold. 
By definition of a regular point, the point $m$ admits a neighborhood ${\mathcal U}$ on which the distribution $E:=\pi^\# (T^*M)$ 
is of constant rank, rank which is an even number, say $2s$.
Since the bivector field $\pi$ is foliated, this distribution is integrable. We can make use of the classical Frobenius theorem
for distributions of constant rank, and deduce that there is a regular foliation ${\mathcal E}$ such that
 $ {\mathcal E}_{m'}=E_{m'} $ for all $m' \in {\mathcal U}$. Of course, $\pi$ is tangent to ${\mathcal E} $
 and its restriction to every leaf of ${\mathcal E}$ is a bivector field of maximal rank.
 It can therefore be inverted to yield, on each leaf of $ {\mathcal E}$, a $2$-form $\omega_{\mathcal E}$.
 In turn, these $2$-form can be extended to a global $2$-forms $\omega$ defined on ${\mathcal U} $:
 For instance one can choose an arbitrary Riemannian metric (which exists when $M$ is oriented) and set $\omega$ to be zero on $ E^\perp$
 while it coincides with $\omega_{\mathcal E}$ on $E$.
We claim that $(\pi,\phi:=\diff \omega) $ is twisted Poisson. It is clear that $\phi$ is a closed $3$-form.
The second condition follows from \cite{CattaneoAlbertoPing}.
\end{proof}

\begin{corollary}
A foliated bivector field of rank two is Poisson.
\end{corollary}
\begin{proof}
 Let $M$ a variety and $\pi$ a foliated bivector field of rank 2. Every regular point $m$ 
 admits, by Frobenius theorem, a neighborhood equipped with   a system of local coordinates $(x_1,\dots ,x_n)$ such that:

  $Im (\pi^\#) = \can{\frac{\partial}{\partial x_1},\frac{\partial}{\partial x_2}}$.
 Consequently $\pi = \varphi(x) \frac{\partial}{\partial x_1} \wedge \frac{\partial}{\partial x_2}$ 
for some smooth function $\varphi$ and a direct computation gives
 $\lb{\pi,\pi}=0$. Regular points being dense in $M$, the relation $\lb{\pi,\pi}=0$ holds at all points.
\end{proof}

There are classes of foliated bivector fields which are a priori not twisted Poisson at non-regular points.
We introduce one of them:

Let $M$ be a smooth manifold. A bivector field $\pi$ is called \emph{a conformally Poisson structure} \cite{BalseiroPaula} if it can be written as $\pi=\phi \pi'$
with $\phi \in \C^\infty(M)$ and $\pi'$ a Poisson structure.
We prove the following:
\begin{proposition}\label{prop:conformally}
Let $\pi$ be a foliated bivector field on a manifold $M$.
Then for every function $\phi \in C^\infty (M)$, the bivector field $\phi \pi$ is
foliated. In particular, conformally Poisson structures are foliated bivector fields.
\end{proposition}
\begin{proof}
Let $\pi'$ be a foliated bivector field, $\phi$ a smooth function on $M$ and $\pi =\phi \pi'$. 
For $F$ a function, we denote by $\X_F$ and $\X_F^{'}$ its Hamiltonian vector field with respect to $\pi$ and $\pi'$ respectively.
The  relation $\X_F=\phi \X_F^{'}$ holds.
Since $\pi'$ is foliated, lemma \ref{lem:immediate}
 implies that there exists a $1$-form $\alpha_{F,G}'$ such that
 $$ \lb{\X_F^{'},\X_G^{'}}= (\pi')^\#(\alpha_{F,G}'),$$
relation that we use to go from the second to the third line of the following computation:
\begin{eqnarray*}
\lb{\X_F,\X_G}
  &=   &\lb{\phi\X_F^{'},\phi\X_G^{'}}\\
   &=   &\phi^{2}\lb{\X_F^{'},\X_G^{'}}+\phi\X_F^{'}[\phi]\X_G^{'}-\phi\X_G^{'}[\phi]\X_F^{'}\\
    &=   &\phi^{2}(\pi')^\#(\alpha_{F,G}') +\phi\X_F^{'}[\phi]\X_G^{'}-\phi\X_G^{'}[\phi]\X_F^{'}\\
      &=   & \phi(\pi')^\#(\alpha_{F,G}')+\X_F^{'}[\phi]\X_G-\X_G^{'}[\phi]\X_F\\
     &=   & \pi^{\sharp}(\phi \alpha_{F,G}'+\X_F^{'}[\phi]\diff G-\X_G^{'}[\phi]\diff F)\\  
      &=   & \pi^\sharp (\alpha_{F,G})  
 \end{eqnarray*}
with $\alpha_{F,G}=\phi \alpha_{F,G}'+\X_F^{'}[\phi]\diff G-\X_G^{'}[\phi]\diff F$.
The result then follows from lemma \ref{lem:immediate}. 
\end{proof}

\begin{example}
We give an example of a bivector field that is not foliated.
We suppose that $M=\R^{3}$ with coordinates $(x,y,z)$ and we set: 
\begin{equation}
\pi = x \pp{}{x} \wedge \pp{}{y} + \pp{}{z} \wedge \pp{}{y}+\pp{}{x}\wedge \pp{}{z}.
\end{equation}
We compute: 
$$\pi^{\sharp }(\diff x)=-x\pp{}{y}-\pp{}{z}~;~ \pi^{\sharp }(\diff y)=x\pp{}{x}+\pp{}{z}~;~ \pi^{\sharp }(dz)=\pp{}{x}-\pp{}{y}.$$
Therefore:
$$\lb{\pi^{\sharp }(\diff y),\pi^{\sharp }(\diff z)}=-\pp{}{x}.$$
There  is no point in $\R^{3}$ where $\pp{}{x}$ is in  the image of $\pi^{\sharp }$. In fact $ x\pp{}{x}+\pp{}{z}$ and $-x\pp{}{x}-\pp{}{y}$ generate the image of $\pi^{\sharp }$.
So if $\pp{}{x}$ in the image, then so are $\pp{}{y}$ and $\pp{}{z}$ hence $\pi^{\sharp }(T^{*}M)=TM$ which is impossible because the image of $\pi^ \#$ must be of even dimension at all point. 
\end{example}

\subsection{General theory of foliated bivector fields}

 It is well-known \cite{Weinstein} that every Poisson structure $\pi$ on a manifold 
 $M$ induces a symplectic foliation, i.e. a foliation whose leaves 
 come with a symplectic structures.
 For foliated bivector fields, there is a similar result, which extends
 a result already well-know for twisted Poisson structures, see 
 \cites{SeveraWeinstein,KosmannSchwarzbachYvette}:
 
 \begin{theorem}\label{theo:localstructureFoliated}
Let $\pi$ be foliated bivector field on a manifold $M$.  
Then 
\begin{enumerate}
 \item There exists a foliation ${\mathcal F}$ on $M$
 such that, for all $m \in M$, the tangent space of the leaf through $m$ 
 is $ \pi_m^\#(T^*M)$
 \item Each leaf of  ${\mathcal F}$ comes with a non-degenerate $2$-form $\omega$ such that for all  $H \in C^\infty(M)$,$\forall m \in M,  \omega(\X_H{_{\vert m}},\bullet )= \diff_mH{_{\vert \Sigma_m}}$ when $\Sigma_m$  is the leaf through $m \in M$
\end{enumerate}
 \end{theorem}
\begin{proof}
The first point follows from the fact that $ \pi^\# (\Omega^1(M))$ 
is a locally finitely generated $C^\infty(M)$-module closed under bracket, so that
the theorem of Sussmann \cite{Sussmann} on integrability of non-regular distribution holds.
By construction, $\pi$ is tangent to each of these leaves and its restriction is of maximal rank. It can therefore 
be inverted to yield a non-degenerate $2$-form.
\end{proof}
Let $E \to M$ and $F \to M $ two vector bundles over the manifold $M$, $P_1$ and $P_2$ two morphisms of vector bundles from $E$ to $F$ over the identity over $M$.
We shall say that the image of $P_1$ is contained in the image of $P_2$, and we write $Im(P_1) \subset Im(P_2) $ if for any local section $\alpha $ of $ E \to M$, there is a section $\beta  $ of $ E \to M$ such that:
 $$  P_1(\alpha ) = P_2 (\beta) .$$

\begin{remark}
\label{rmk:Image}
Note that if $ Im(P_1) \subset Im(P_2)$, then at any point $m \in M$, the image of the linear mapping $ (P_1)_m: E_m \to F_m$ is contained in the image of the linear mapping $(P_2)_m: E_m \to F_m$.
The converse is false as shown by the following example:
 $$\pi_{A}=(x^{2}+y^{2})\pp{}{x} \wedge \pp{}{y} \hbox{ and } \pi_{B}=(x^{2}+y^{2})^{2}\pp{}{x} \wedge \pp{}{y}.$$
It is clear that $Im(\pi^{\sharp}_{A})_m = Im(\pi^{\sharp}_{B})_m$ at every point $m \in \R^{2}$. However
for $\alpha = \diff x$, the vector field
$$\pi^{\sharp }_{A}(\alpha) = (x^{2}+y^{2})\pp{}{y} $$
can not be written as $\pi_B^\# (\beta)$, since a $1$-form $\beta = F(x,y) \diff x+G(x,y) \diff y$ such that $\pi^{\sharp}_{B}(\beta) = \pi^{\sharp}_{A}(\alpha)$ should satisfy:
$$\pi^{\sharp }_{B}(F(x,y)\diff x+G(x,y)\diff y)=(x^{2}+y^{2})^{2}F(x,y)\pp{}{y}+(x^{2}+y^{2})G(x,y)\pp{}{x}. $$
This imposes, $G(x,y)=0$ and $(x^{2}+y^{2})^{2}F(x,y)=(x^{2}+y^{2})$, and $F(x,y)= \frac{1}{x^{2}+y^{2}}$ for all  $(x,y) \neq (0,0)$.
But $F(x,y)= \frac{1}{x^{2}+y^{2}}$ can not extended by continuity at the point $(0,0)$ so such a $1$-form $\beta$
does not exist. 
\end{remark}

Together with lemma \ref{lem:immediate}, the next criterium shall be useful in the sequel:

\begin{proposition}\label{prop:foliated=image}
Let $\pi$ a bivector field  on a manifold $M$. Then the following points are equivalent: 
\begin{enumerate}
\item[(i)] $\pi$ is foliated,
\item[(ii)] for any function $H \in C^\infty(M)$, we have $Im(({\mathcal L}_{\X_H} \pi)^\#) \subset Im(\pi^\#). $
\end{enumerate}
\end{proposition}
\begin{proof}
Let $H  \in C^\infty(M)$ and 
 $\alpha \in \Omega^1(M)$. 
 \begin{eqnarray*}
\lb{\X_H , \pi^\# ( \alpha ) } &=& {\mathcal L}_{\X_H}  ( \pi^\# (\alpha ) ) \\
 & = & ({\mathcal L}_{\X_H}  \pi)^\# (\alpha)+ \pi^\# ({\mathcal L}_{\X_H} \alpha ). 
 \end{eqnarray*}  
  If $\pi$ is foliated, then there exist $\omega \in \Omega^1(M)$ such that 
$\lb{ \pi^\# (\diff H) , \pi^\# (\alpha) } = \pi^\# (\omega )$,
 hence 
$$ ({\mathcal L}_{\X_H} \pi)^\# (\alpha) =   (\pi^\#)(\omega) - \pi^\# ({\mathcal L}_{\X_H} \alpha)
= \pi^\# (\omega - {\mathcal L}_{\X_H} \alpha), $$
 so that $ Im ({\mathcal L}_{\X_H} \pi)^\# \subset Im (\pi^\#)$ and \emph{(i)} implies \emph{(ii)}.
 Conversely, if the inclusion $Im(({\mathcal L}_{\X_H} \pi)^\#) \subset Im(\pi^\#)$ holds, then 
  for any function $F \in C^\infty (M)$ there is a $1$-form $\beta_F$ such that $({\mathcal L}_{\X_H} \pi)^\#\diff F=\pi^\#(\beta_F)$.
  The previous relation applied to $\alpha = \diff F$,
   \begin{eqnarray*}
\lb{\X_H , \X_F} &= &\lb{ \X_H ,\pi^\# (\diff F) } \\
 & = & ({\mathcal L}_{\X_H} \pi )^\# \diff F + \pi^\# ({\mathcal L}_{\X_H} \diff F) \\
 & = &  \pi^\# (\alpha_{F,H}) \hbox{ with } \alpha_{F,H} := \beta_F + {\mathcal L}_{\X_H} \diff F.
 \end{eqnarray*}  
Since $F$ and $H$ are arbitrary smooth functions on $M$, lemma \ref{lem:immediate} implies that $\pi$ is a foliated bivector field
and \emph{(ii)} implies \emph{(i)}. This completes the proof.
\end{proof}

\begin{proposition}\label{prop:11tensor}
Let $(M,\pi$) be a foliated bivector field. Then for every function $H \in {\C}^\infty(M)$, there exists a $1$-$1$ tensor $C_H: TM \to TM $ such that 
 $  ({\mathcal L}_{\X_H} \pi)^\# = C_H \circ \pi^\#$.
\end{proposition}
The proof uses the following lemma:
\begin{lemma}\label{lem:generalitesMorphismes}
Let $E \to M$ and $F \to M $ be a pair of vector bundles over $M$.
Let $P_1$ and $P_2$ two morphisms of vector bundles over the identity of the manifold $M$. 
If $Im(P_1) \subset Im(P_2) $, then there exists $ D$ such that $ P_1 =P_2 \circ D$.
\end{lemma}
\begin{proof}
We begin by showing the result locally.
Let $m \in M$ be a point, and ${\mathcal U} $ a neighborhood of this point on which there is a local trivialization $(e_1, \dots,e_d)$ of $E$, with $d$ the rank of $E$.
Then for all $i=1, \dots, d$ there exists a section $\beta_i$ of $E$ such that $P_1(e_i)=P_2(\beta_i)$.
As the sections $(e_1, \dots,e_d)$ form a local basis of $E$, we can rewrite the sections $\beta_i$ in this basis, thus obtaining: $$ \beta_i = \sum_{j=1}^d H_i^j e_j$$
 therefore
  $$P_1 (e_i) = P_2 \circ D_{\mathcal U} (e_i) $$
where $D_{\mathcal U} $ is the matrix $(H_i^j)_{i,j=1,\dots,d} $, 
matrix corresponding to a morphism of vector bundles over the identity of ${\mathcal U}$ still denoted by $D_{\mathcal U}$.
The identity  $ P_2 \circ D_{\mathcal U} =P_1$ holds on ${\mathcal U} $ by construction.
Now, choose an open cover $({\mathcal U}_k)_{k \in K}$ of $M$ by open sets associated a partition of unity $(\varphi_k)_{k \in K}$, such that $E$ is trivializable
on each open subset ${\mathcal U}_k$. For each indice $k \in K$, there exists
 a morphism of vector bundles $D_k : E \to E$, over the identity of ${\mathcal U}_k$, such that $ P_2 \circ D_k=P_1$. The morphism of vector bundle $D: E \to E$, over the identity of $M$, given by
$ D:=\sum_{k \in K} \varphi_k D_k $ satisfies $P_1 = P_2 \circ D $ by construction.
The lemma follows.
\end{proof}

We can now prove the proposition \ref{prop:11tensor}.

\begin{proof}
By proposition \ref{prop:foliated=image}, $\pi$ is foliated if and only if 
 $  Im\left({\mathcal L}_{\X_H} \pi \right)^\# \subset  Im(\pi^\#) $ for an arbitrary smooth function $H$.
By lemma \ref{lem:generalitesMorphismes}, there exists, for all $H \in C^\infty(M)$,
a morphism $D_H : T^*M \to T^*M$ of vector bundles over the identity of $M$ such that
\begin{equation}\label{eq:CindiceH} 
 ({\mathcal L}_{\X_H} \pi)^\# =  \pi^\# \circ D_H.
\end{equation}
Both sides of (\ref{eq:CindiceH}) are vector bundle morphisms from $T^*M $ to $TM$, so their dual maps 
also are  vector bundle morphisms from $T^*M $ to $(T^*M)^*=TM$ (upon identifying the bidual with the dual).
Recall that for every bivector field $\pi'$,  the morphism $(\pi ')^\#: T^*M \to TM $ is antisymmetric, i.e.
 $  ((\pi')^\#)^* = - (\pi')^\# $.
Taking the dual of (\ref{eq:CindiceH}) therefore yields
$ ({\mathcal L}_{\X_H} \pi)^\# =  D_H^* \circ \pi^\# $ where $ D_H^* : TM \to TM$ is the dual of $D_H$. Note that 
$D_H^*$ is a $1$-$1$ tensor, so that $C_H =  D^*$ satisfies the requirements of proposition \ref{prop:11tensor}.
\end{proof}

\subsection{Weakly foliated bivector field}

As we saw in the course of remark \ref{rmk:Image}, a vector field $X$ on a manifold $M$ equipped 
with a bivector field $\pi$ may satisfy $X_m \in Im(\pi_m^\#)$ for every point $ m \in M$
without being of the form $X=\pi^\#(\beta)$ for some $1$-form $\beta \in \Omega^1(M)$.
This difference is at origin of the following definition.

\begin{definition}
A bivector field $\pi$ on a manifold $M$ said to be weakly foliated when for any pair  $\alpha,\beta \in \Omega ^{1}(M)$, and any point $m \in M$, we have:
$ \lb{\pi^\# (\alpha) , \pi^\# (\beta)}{_{\vert m}}  \in Im(\pi^\#_m)$.
\end{definition}

\begin{example}
Foliated bivector fields are weakly foliated.
\end{example}

\begin{example}\label{example:weaktwisted}
Proposition \ref{prop:twistedmeansfoliated} generalizes as follows.
If at every point $m \in M$, we have
 $[\pi,\pi]_m \in Im( \wedge^3 \pi_m^\#)$ (which is always the case for a twisted Poisson structure),
 then $\pi$ is weakly foliated.

 The proof goes as follows.
 First, lemma \ref{lem:immediate} can be easily generalized: 
 a bivector field is weakly foliated if and only if $\lb{\X_F,\X_G}{_{\vert m}} \in Im(\pi^\#_m)  $ for all $ F,G \in \F(M)$ and $m \in M$. 
If $[\pi,\pi]_m \in  \wedge^3 \pi_m^\#(\omega)$ for some $\omega \in \wedge^3(T_m^* M)$, then,
for every functions $F,G$ and $1$-form $\alpha$:
\begin{eqnarray*} \can{ \lb{\X_F,\X_G},\alpha }{_{\vert m}} &=& \can{\X_{\pb{F,G}}{_{\vert m}},\alpha }+ \can{ \lb{\pi,\pi}_m, \diff_m F \wedge \diff G_m \wedge \alpha_m  }  \\
  & = & \can{\X_{\pb{F,G}}{_{\vert m}},\alpha }+\can{\wedge^3 \pi^\# (\omega ),\diff_m F \wedge \diff_m G \wedge \alpha,  } \\
  & = & \can{\pi^\#_m ( \diff_m \pb{F,G} ) ,\alpha }+\can{\pi^\#_m (\diff_m F) \wedge \pi^\#_m (\diff_m G) \wedge \pi^\#(\alpha), \omega } \\
  & = & \can{\pi^\#_m ( \diff_m \pb{F,G} ) ,\alpha }+\can{\pi^\#(\alpha), \imath_{(\X_F )_m \wedge (\X_G )_m} \omega } \\
  & = & \can{\pi^\#_m ( \diff_m \pb{F,G} ) ,\alpha }-\can{\alpha, \pi^\#\left(\imath_{(\X_F)_m \wedge (\X_G )_m } \omega\right) } \\
  & = & \can{\pi^\#_m( \diff_m \pb{F,G} -\imath_{( \X_F )_m \wedge (\X_G )_m } \omega ) ,\alpha }.
\end{eqnarray*}
Above, we have used, to go from the first to the second line the relation:
 $$ \lb{\X_F,\X_G}-\X_{\pb{F,G}}=\imath_{\diff F \wedge \diff G} \lb{\pi,\pi},$$
 while the rest of the computations are pure multilinear algebra.
The previous relations being valid for all $\alpha$, it implies  
$$\lb{\X_F,\X_G}{_{\vert m}} = \pi^\#_m( \diff_m \pb{F,G} -\imath_{(\X_F)_m \wedge (\X_G)_m } \omega),$$
which proves the claim.
\end{example}

We leave it to the reader to adapt the proof  of proposition  \ref{prop:foliated=image} to yield:
\begin{proposition}\label{foliated=imageweakly}
A bivector field is weakly foliated if and only if
for every function $F \in C^\infty(M)$, and every point $m \in M$, the inclusion 
$ Im ({\mathcal L}_{\X_F} \pi)_m^\# \subset Im (\pi_m^\#)$ holds.
\end{proposition}
\begin{example}\label{ex:WeaklyNotStrongly}
We give two examples of weakly foliated bivector fields which are not foliated. Let $N$ a manifold. Consider the bivector field on $M:=N \times \R^{2}$ 
defined by: 
\begin{equation}
\pi = X \wedge \pp{}{u} +Y \wedge  \pp{}{v}
\end{equation}
where $u,v$ are the coordinates on $\R^{2}$ and $X,Y$ are two vector fields on $M$.
By construction:
\begin{equation}\label{eq:pisharpisXisY}
\pi^{\sharp }(\diff u)=X \hbox{ and } \pi^{\sharp }(\diff v)=Y .
\end{equation}
We now choose $N=\R^{2}$ endowed with the canonical coordinates $(x,y)$ and we set:
$$X=  x^i\pp{}{x} \hbox{ and } Y=(x^{2}+y^{2})\pp{}{y}$$
where $i$ is either $0$ or $1$, so that
\begin{equation}\label{eq:explicitNotFoliated}
\pi =  x^i\pp{}{x}  \wedge \pp{}{u} +(x^{2}+y^{2})\pp{}{y} \wedge  \pp{}{v}.
\end{equation}
The relation
$$\lb{X,Y} = \lb{x^i\pp{}{x},(x^{2}+y^{2})\pp{}{y}} = 2x^{i+1}\pp{}{y} =  \frac{ 2x^{i+1}}{ x^{2}+y^{2}}Y, $$
together with (\ref{eq:pisharpisXisY}), implies that, for all $(x,y) \neq (0,0)$ :
\begin{equation}\label{eq:dudvbeta}
\lb{\pi^{\sharp }(\diff u),\pi^{\sharp }(\diff v)}=\pi^{\sharp }(\beta) \hbox{ with }  \beta = \frac{ 2x^{i+1}}{ x^{2}+y^{2}}\diff v.
\end{equation}
If $x \neq 0$, the bivector field $\pi$ is invertible at  $(x,y,u,v) \in M$,
so that $\beta=\frac{ 2x^{i+1}}{ x^{2}+y^{2}} \diff v$ is the unique covector that satisfies (\ref{eq:dudvbeta}).
But $\frac{ 2x^{i+1}}{ x^{2}+y^{2}}dv$ can not be extended by continuity at a point of the form  $(0,0,u,v) \in M$ and consequently
the bivector field $\pi$  is not foliated.
We claim, however, that it is weakly foliated. This can be deduced from the relation
 \begin{eqnarray*} \lb{\pi,\pi}& = & 2 \lb{X,Y} \wedge \pp{}{u} \wedge  \pp{}{v} \\
  & = &  2  \frac{ 2x^{i+1}}{ x^{2}+y^{2}} Y\wedge \pp{}{u} \wedge  \pp{}{v} \\
  & = &  4  \frac{ x}{ (x^{2}+y^{2})^2}  \pi^\# (\diff v) \wedge  \pi^\# (\diff x) \wedge \pi^\# (\diff y ) \\
  & = & \wedge^3 \pi^\# (\omega )
 \end{eqnarray*}
with $\omega =  4  \frac{x}{ (x^{2}+y^{2})^2} \diff v \wedge  \diff x \wedge \diff y $.
 The criterion  proposed in example \ref{example:weaktwisted} is satisfied if $(x,y) \neq (0,0)$, implying that $ \pi$ is weakly foliated
 at these points. Now, for $i=1$, at a point of the form $m=(0,0,u,v) \in M$,
 $\pi_m =0$, which implies that $  \pi $ is weakly foliated on ${\mathbb R}^4$.
 For $i=0$,  at a point of the form $m=(0,0,u,v) \in M$, the image of $\pi^\#$ is $\pp{}{x},\pp{}{u} $,
 which is an integrable distribution, hence the result holds true at these points also. 
\end{example}

\section{Poisson structures and foliated bivectors: Lagrangian theory}
\label{sec:LagrangianTheory}

Out of a bivector field $\pi$ and a function $H$ on a manifold $M$, We define a Lagrangian on the set of all paths valued in $T^*M $.
Stationary points of this Lagrangian are shown to depend on the properties of the Jacobiator of $\pi$.
In section \ref{sec:defnot}, we define subsets of the set of paths, in the continuation of \cite{CattaneoFelder}.
In section \ref{sec:defFunct}, the functional is introduced and its differential defined and computed.
In section \ref{sec:mainresults}, we state and prove the main theorem of the paper and we give in section \ref{sec:counterexamples} two counterexamples
that prevent us from getting a better result.

\subsection{Definitions and notations}
\label{sec:defnot}

Let $(M,\pi)$ be a manifold equipped with a bivector field $\pi$ (which is not assumed to be Poisson). 
We consider the set $\tilde{P}(T^* M)$ of smooth paths from $I=[0,1]$ to $T^*M$.

Recall \cite{CrainicFernandes} that a \emph{cotangent paths for $\pi$} is a path $\alpha \in \tilde{P}(T^* M)$  which satisfies:
  \begin{equation} \label{eq:def_chemins_cotangents}
 \pi^\# (\alpha(t)) = \frac{\diff x(t)}{\diff t}   
  \end{equation}
where $x = p \circ \alpha$, called  \emph{base path}, is the projection onto $M$ of the path $\alpha(t)$.
We denote by ${P}(T^* M)$ the set of cotangent paths for $\pi$, following the convention of \cites{CrainicFernandes,FernandesCranicMarius}.

We need to define another class of paths that contains ${P}(T^* M)$.
We start by introducing a class of paths in $TM$ associated to a vector field.
Given a vector field $X$  of flow $ \varphi_t$ on a variety $M$, we call  \emph{ tangent integral curve for $X$} a smooth path 
$ \beta $ from $ I=[0,1]$ to $TM$ 
 that verifies:
   \begin{equation}\label{eq:tangentintcurve} 
b(t) = T\varphi_{t} (b (0)).
 \end{equation}
Equation (\ref{eq:tangentintcurve}) implies that the base path of $ b(t)$ is an integral curve of $X$.

\begin{lemma}\label{lem:int_tang}
Let $M$ be a manifold equipped with an arbitrary connection $\nabla$ on $TM$ of torsion $\Tor^\nabla$.
For all $X$ vector field on $M$, a path $b(t) : I=[0,1] \to TM$ is a tangent integral curve of $X$ if and only if:
\begin{enumerate}
\item  The base path of $b(t)$ is on integrable curve of $X$,
\item   $\nabla_X b = \nabla_b X + \Tor^\nabla (X,b)$.
\end{enumerate}
(Note that condition $1$ allows to make sense of the quantity $\nabla_X b$.)
\end{lemma}
\begin{proof}
The result is clearly true on $\R^n$ equipped with the canonical connection, having in mind the fact that the base path of 
a tangent integral curve of $X$ is an integral curve of the vector field $X$. 
It is then easy to see that the quantity $\nabla_X b - \nabla_b X + \Tor^\nabla (X,b)$ does not depend on 
the connection $\nabla$, which completes the proof.
\end{proof}

Returning now to the general case, without assuming that $M$ is included in an open of ${\mathbb R}^n$. 
We shall make use of the following notion.
\begin{definition}
Let $M$ be a manifold.
equipped with a bivector field $\pi$ and a Hamiltonian function $H \in C^\infty(M)$,
we call  \emph{quasi-cotangent path} for $(\pi,H)$ a path $ \alpha \in \tilde{P}(T^* M)$,
of base path $x : I=[0,1] \to M$,
such that $c(t) :=  \pi^\#_{x(t)} (\alpha(t)) - \frac{\diff x(t)}{\diff t}   $  
is an tangent integral curve of the Hamiltonian vector field $\X_H$. 
\end{definition}

\begin{example}
Of course, cotangent paths for $\pi$ are for every Hamiltonian function $H \in C^\infty(M)$ quasi-cotangent paths for $(\pi,H)$.
In particular,  the path
 $$ t \mapsto \diff_{x(t)} H ,$$
where $x(t)$ is an integral curve of $X_H$, is both a quasi-cotangent path for $(\pi,H)$ and a cotangent path for $\pi$.
\end{example}

The following lemma follows from lemma \ref{lem:int_tang}.

\begin{lemma}\label{lem:quasi_hami}
Let $M$ be a manifold equipped with a bivector field $\pi$
and an arbitrary connection $\nabla$.
 For every function $H \in C^\infty(M)$, define a $1$-$1$ tensor $K^{H} $ by 
\begin{equation} \label{eq:defK} K^H_x(u):=\left. \nabla_u \X_H \right|_x+Tor^\nabla (\X_H ,u)\end{equation}
for all $x\in M$ and all tangent vector $u \in T_x M$.
A path $\alpha \in \tilde{P}(T^*M)$ 
is a quasi-cotangent path for $(\pi,H)$
if and only if its base path $x(t)$ follows the flow of $\X_H$ and
  $$ \nabla_{\frac{\diff x}{\diff t}} c =K^H (c)$$
  where $c(t) = \pi^ \#_{x(t)} (\alpha(t))  - \frac{\diff x(t)}{\diff t}.$
\end{lemma}

\subsection{A functional on the space of paths valued in the cotangent bundle}
\label{sec:defFunct}

Let $M$ be a manifold equipped with a bivector field $\pi$.
To any function $H \in \FF(M)$, we associate ${\mathbb R}$-valued function on $\tilde{P}(T^*M)$ 
by:
  \begin{equation}
   \label{eq:defL}
 {\mathcal L}^H (\alpha) = \int_0^1 \left\langle \left. {\mathcal X}_H \right|_{x(t)}  -  \frac{\diff x(t)}{\diff t} , \alpha(t) \right\rangle \diff t ,
  \end{equation}
 for all path $\alpha \in \tilde{P}(T^*M)$ with base path $x$.
 
This functional is related to several functionals that appear in the literature, but is however different,
see section \ref{sec:LagrangianTheory}.
 
We would like to study the stationary points of this functional. For this purpose, we first need to say it which sense it is differentiable.
Although it should be possible to deal with $C^1$-paths instead of smooth paths \cite{CattaneoFelder}, therefore placing ourself within the context
of infinite dimensional geometry, we prefer to work with smooth paths, and to speak of G\^ateau differentiability.

Let $N$ be a manifold. Let $\gamma \in \tilde{P}(N)$ be a smooth path. We denote by  $T_\gamma (\tilde{P}(N)) $ and call 
\emph{tangent space of $ \tilde{P}(N)$ at the point $\gamma$} the vector space of all smooth maps $e$ from $I=[0,1] $ to $T(N)$
with base bath $\gamma$ (the base path here meaning the image of $e$ through the projection $TN \to N$.
For every \emph{deformation $(\epsilon,t) \to \gamma_\epsilon(t)$ of $\alpha$}, i.e. every smooth map 
$$\begin{array}{rcl} 
[0,1] \times ]-u, +u[ &\mapsto& N \\
(t,\epsilon) & \to & \gamma_\epsilon (t) 
  \end{array}
$$
such that $\gamma_0(t) = \gamma(t)$, notice that $ t \to \left. \frac{\diff \gamma_\epsilon (t)}{ \diff \epsilon}\right|_{\epsilon=0}$
is an element in the tangent space of  $ \tilde{P}(N)$ at the point $\gamma$ (by construction, it is a map from $I=[0,1]$ to $TN$ above the path $\gamma: I \to N$).
We say that a function ${\mathcal L}$ from $  \tilde{P}(N)$ to ${\mathbb R}$ is \emph{differentiable} at a point $\gamma$ 
if there exists a linear form $\diff_\gamma {\mathcal L}$ on the tangent space $T_\gamma (\tilde{P}(N))$ such that
  $$\diff_\gamma {\mathcal L} \left( \left. \frac{\diff \gamma_\epsilon (t)}{ \diff \epsilon}\right|_{\epsilon=0}  \right)  = \left.\frac{\diff  {\mathcal L} 
  \left(\gamma_\epsilon \right)}{ \diff \epsilon}\right|_{\epsilon=0}$$
for all deformation $(\epsilon,t)\to  \gamma_\epsilon(t) $ of $\gamma$. We call $\diff_\gamma {\mathcal L} $ the \emph{differential of $ {\mathcal L}$ at $\gamma$}
and we say that $\gamma$ is a \emph{stationary point} when this differential is zero.
We say that $\gamma \in S \subset \tilde{P}(N)$ is a \emph{stationary point when restricted to some subset  $S$ of $\tilde{P}(N)$}
when the differential vanishes on every tangent vector of the form $ \left. \frac{\diff \gamma_\epsilon (t)}{ \diff \epsilon}\right|_{\epsilon=0}$
with  $(\epsilon,t)\to  \gamma_\epsilon(t) $ a deformation of $ \gamma$ valued in $S$, i.e. such that the path $\gamma_\epsilon$
is in $S$ for all $\epsilon \in ]-u,u[$. (In practice, our subsets shall of course be, at least morally, infinite dimensional submanifolds
of $\tilde{P}(N)$).

To express in an explicit manner the differential of ${\mathcal L}^H$, we 
choose an arbitrary connection $ \nabla$ on $M$. This allows to identifying
the tangent space  $T_\alpha(T^*M) $ with $T_xM \oplus T_x^*M$ where $x=p(\alpha )$ is the base point of $\alpha \in T^*M$.
Upon choosing a connection $\nabla$, we have the natural identification:
\begin{equation}\label{eq:isomNabla}
 T_\alpha \tilde{P}(T^*M) \simeq_\nabla \Gamma(x^*(TM \oplus T^*M))
\end{equation}
More precisely, one can identify an element in $T_\alpha \tilde{P}(T^*M)$ with a pair $(\gamma_0,\delta_0) $ of smooth maps from $I=[0,1]$ to $TM$ and $ T^*M$ respectively
such that for all $t \in I$, $\gamma_0(t) $  and $\delta_0 $ belongs to $ T_{x(t)}M$ and $ T^*_{x(t)}M$ respectively.
  
\begin{proposition} \label{prop:calcul_diff}
Let $M$ be a manifold equipped with a bivector field $\pi$. 
For every function $H \in C^\infty(M)$, the functional ${\mathcal L}^H$ is differentiable (in the sense above) at all point  $\alpha \in \tilde{P}(T^*M)$.

Moreover, upon choosing a connection $\nabla$ with torsion $\Tor^\nabla$, and identify an element in $T_\alpha \tilde{P}(T^*M)$ with a pair $(\gamma_0,\delta_0) $,
as in (\ref{eq:isomNabla}), the differential of $ {\mathcal L}^H $ at  the point $ \alpha \in \tilde{P}(T^*M) $ in the direction 
 of $(\gamma_0,\delta_0) \in T_\alpha \tilde{P}(T^*M)$
 is then given by:
 \begin{eqnarray*}
 \diff_{\alpha} {\mathcal L}^H (\gamma_0,\delta_0)
  &=   &\int_0^1 \left\langle \delta_0(t),   \left. {\mathcal X}_H \right|_{x(t)} - \frac{\diff x(t)}{\diff t} \right\rangle \diff t
\\&&+\int_0^1 \left\langle \alpha ,\nabla_{\gamma_0} \X_H + \Tor^\nabla \left(\frac{\diff x (t)}{\diff t} ,\gamma_0 \right) 
+ \nabla_{\frac{\diff x (t)}{\diff t}} \gamma_0  \right\rangle \\
 &=   &\int_0^1 \left\langle \delta_0(t),   \left. {\mathcal X}_H \right|_{x(t)} - \frac{\diff x(t)}{\diff t} \right\rangle \diff t
 \\&&+\int_0^1 \left\langle \alpha ,\nabla_{\gamma_0} \X_H + \Tor^\nabla\left(\frac{\diff x (t)}{\diff t} ,\gamma_0   \right)\right\rangle \\
\\&&-\int_0^1 \left\langle  \nabla_{\frac{\diff x (t)}{\diff t}} \alpha , \gamma_0   \right\rangle \\&&
 + \left\langle \alpha(1),\gamma(1) \right\rangle - \left\langle \alpha(0),\gamma(0) \right\rangle. 
 \end{eqnarray*}
In the previous, we used the same notation $\nabla$ for a connection on $TM$ and its induced connection on $T^*M$.
\end{proposition} 
\begin{proof}
To go from the first to the second equality of the display of proposition \ref{prop:calcul_diff},
we use integration by parts, which takes, in this context, the following form:
for a given path $x$ on $M$, and arbitrary $\beta : I \to  T^*M, \xi : \to TM$ such that $\beta(t) \in T^*_{x(t)} M,\xi(t) \in T_{x(t)}M$ for all $t \in I$, the 
following relation holds:
$$  \int_0^1  \left\langle  \nabla_{\frac{\diff x (t)}{\diff t}} \beta , \xi   \right\rangle dt+ \int_0^1  \left\langle  \beta, \nabla_{\frac{\diff x (t)}{\diff t}} \xi \right\rangle dt
= \langle  \beta (1), \xi (1)   \rangle - \langle  \beta (0), \xi (0)  \rangle  .$$
It suffices therefore to prove the first formulation of the differential of ${\mathcal L}^H $.
Let  $(\epsilon,t) \to \alpha_\epsilon(t)$ be a deformation of $\alpha$. 
By construction, $ \left. \frac{\diff \alpha_\epsilon (t)}{ \diff \epsilon}\right|_{\epsilon=0} $ is a path valued in $T(T^*M)$
above $\alpha$, which in view of (\ref{eq:isomNabla}) can be identified with an element $( \gamma_0,\delta_0)$ in 
$\Gamma(x^*(TM \oplus T^*M))$. 
According to the definition of ${\mathcal L}^H$ in equation (\ref{eq:defL}), and Fubini theorem, we have:
$$ \left. \frac{\diff  {\mathcal L}^H \left(\alpha_\epsilon \right)}{ \diff \epsilon} \right|_{\epsilon=0} = 
 \int_0^1 \left. \frac{\diff}{\diff \epsilon}\right|_{\epsilon=0} \left\langle \left. {\mathcal X}_H \right|_{x_\epsilon(t)}  -
 \frac{\diff x_\epsilon(t)}{\diff t} , \alpha_\epsilon(t) \right\rangle \diff t$$
 where $t \to x_\epsilon$(t) is for every $\epsilon  \in ]-u,u[ $  the base path of the path $ t\to \alpha_\epsilon(t)$.
 To complete the proof of the proposition, it suffices to establish the following identity for all $t \in I$:
 %We wish to establish for all $t \in [0,1]$ the identity:
$$ \begin{array}{ll}
&\left. \frac{\diff}{\diff \epsilon}\right|_{\epsilon=0} \left\langle \left. {\mathcal X}_H \right|_{x_\epsilon(t)}  -
 \frac{\diff x_\epsilon(t)}{\diff t} , \alpha_\epsilon(t) \right\rangle  \\ 
  =&\left\langle \delta_0(t),   \left. {\mathcal X}_H \right|_{x(t)} - \frac{\diff x(t)}{\diff t} \right\rangle \\
+ &\left\langle \alpha ,\nabla_{\gamma_0} \X_H + \nabla_{\frac{\diff x (t)}{\diff t}} \gamma_0+ \Tor^\nabla\left(\frac{\diff x (t)}{\diff t} ,\gamma_0 \right)  \right\rangle .
 \end{array}$$
This result is obvious when $M={\mathbb R}^n$ is equipped with the canonical connection $\nabla_0 $, 
case in which is just amounts to a simple computation of differential:
\begin{equation}\label{eq:todifferentiate} \begin{array}{ll}
&\left. \frac{\diff}{\diff \epsilon}\right|_{\epsilon=0} \left\langle \left. {\mathcal X}_H \right|_{x_\epsilon(t)}  -
 \frac{\diff x_\epsilon(t)}{\diff t} , \alpha_\epsilon(t) \right\rangle  \\ 
  =&\left\langle \delta_0(t),   \left. {\mathcal X}_H \right|_{x(t)} - \frac{\diff x(t)}{\diff t} \right\rangle \\
+ &\left\langle \alpha ,\diff \X_H(\gamma_0) + \frac{\diff  \gamma_0 (t)}{\diff t}  \right\rangle .
 \end{array}\end{equation}
which is the desired quantity since for the canonical connection, the torsion vanishes, the identities 
$\frac{\diff  \gamma_0 (t)}{\diff t} =(\nabla^0)_{\frac{\diff x (t)}{\diff t}} \gamma_0$ and
 $\diff \X_H(\gamma_0)  = (\nabla^0)_{\gamma_0}(\X_H)$ hold.
 
To prove the general case, it suffices therefore to establish that the quantity that appears on the right hand side of (\ref{eq:todifferentiate}), i.e.
 \begin{equation}\label{eq:doesNotDependConnection} 
 \left\langle \delta_0(t),   \left. {\mathcal X}_H \right|_{x(t)} - \frac{\diff x(t)}{\diff t} \right\rangle 
+ \left\langle \alpha ,\nabla_{\gamma_0} \X_H + \nabla_{\frac{\diff x (t)}{\diff t}} \gamma_0 + \Tor^\nabla\left(\frac{\diff x (t)}{\diff t} ,\gamma_0 \right)  \right\rangle, 
\end{equation}
does not depend on the choice of a connection $\nabla$, since this invariance allows to use local charts to reduce the problem 
to ${\mathbb R}^n$ equipped with the canonical connection where it is shown to be true.

 Let $\nabla'$ be a second connection on $M$.
There exists a $2$-$1$ tensor $A : TM \times TM \to TM$ that we chose to denote by $(u,v) \mapsto A_u (v)$ such that:
 $$ \nabla_u 'v = \nabla_u v + A_u (v)$$
 for all tangent vectors $u,v$ in the same tangent space.
 Recall that the identification of $T_\alpha \tilde{P}(M) $ with pair of paths above the base path $x$ in $TM$ and $T^*M$
depends on the connection.
 More precisely, upon changing $\nabla$ to $\nabla'=\nabla+a$, the isomorphism (\ref{eq:isomNabla}) amounts to the following
 transformation
 $$ (\gamma_0,\delta_0) \to \left(\gamma_0 , \delta_0 + A_{\gamma_0}^*(\alpha) \right).$$
When changing the connections, all the terms adding up to the right hand side of (\ref{eq:doesNotDependConnection})
 are therefore modified and the array below recapitulates how:
$$\begin{array}{lcl}
\left\langle \delta_0(t),   \left. {\mathcal X}_H \right|_{x(t)} - \frac{\diff x(t)}{\diff t} \right\rangle & \hbox{ becomes }  & 
\left\langle \delta_0(t),   \left. {\mathcal X}_H \right|_{x(t)} - \frac{\diff x(t)}{\diff t} \right\rangle   \\ & & + \left\langle A_{\gamma_0}^* (\alpha) , 
\left. {\mathcal X}_H \right|_{x(t)} - \frac{\diff x(t)}{\diff t}  \right\rangle \\
 \left\langle \alpha ,\nabla_{\gamma_0} \X_H   \right\rangle & \hbox{ becomes }  &  \left\langle \alpha ,\nabla_{\gamma_0} \X_H  \right\rangle + 
 \left\langle \alpha ,A_{\gamma_0} \X_H \right\rangle \\
 \left\langle \alpha ,\nabla_{\frac{\diff x (t)}{\diff t}} \gamma_0 \right\rangle & \hbox{ becomes }  &  \left\langle \alpha ,\nabla_{\frac{\diff x (t)}{\diff t}} \gamma_0  
\right\rangle+
\left\langle \alpha ,A_{\frac{\diff x (t)}{\diff t}} \gamma_0  \right\rangle\\
 \left\langle \alpha ,\Tor^\nabla\left(\frac{\diff x (t)}{\diff t} ,\gamma_0 \right)  \right\rangle  & \hbox{ becomes }  &  
 \left\langle \alpha ,\Tor^\nabla\left(\frac{\diff x (t)}{\diff t} ,\gamma_0 \right)  \right\rangle \\ & & + \left\langle \alpha ,A_{\frac{\diff x (t)}{\diff t}} (\gamma_0 ) -
  A_{\gamma_0} \left(  \frac{\diff x (t)}{\diff t}\right) \right\rangle.
\end{array}$$
 Adding up the terms of the first or the third column, one obtains the same quantity, which proves that (\ref{eq:doesNotDependConnection})
 does not depend on the choice of a connection and therefore completes the proof.
\end{proof}

\subsection{Statement of the main results}
 \label{sec:mainresults}
 
We introduce two types of subsets of the set $\tilde{P}(T^* M)$ of all paths from $I=[0,1]$ to $T^{*}M$.

\begin{definition}
Let $M$ be a manifold. Given two points $m,m' \in M$, 
we denote by $\tilde{P}_{m,m'}(T^*M)$ and call  \emph{paths connecting $m$ to $m'$} subset of $\tilde{P}(T^* M)$ which satisfy  $x(0) =m$ and  $x(1) =m'$.
Here, $x = p \circ \alpha $ is on usual mode of paths the base path of $\alpha$.

Assume now $M$ is equipped with a bivector field $\pi$.
We denote by $\hat{P}_{m,m'}(T^*M)$ and call \emph{initially cotangent paths for $\pi$ connecting $m$ to $m'$} the set of paths $\alpha$ connecting $m$ to $m'$ 
such that the equation (\ref{eq:def_chemins_cotangents}) is satisfied for $t=0$, i.e.:
$$\pi^\# (\alpha(0)) = \left. \frac{\diff x(t)}{\diff t}\right|_{t=0}  $$
\end{definition}

Given two points $m,m' \in M$, the function ${\mathcal L}^H (\alpha)$ restricts to a function, still denoted by $ {\mathcal L}^H$,
on each of the previous two subsets. The main purpose of this paper is to prove the following result, interpreting
the stationary point of these restricted functionals:

\begin{theorem}\label{theo:caract_chemins_cotangents}
Let $M$ be a manifold and $\pi $ be a bivector field on $M$.
\begin{enumerate}
\item  If the bivector field $\pi$ is foliated, then
for every function $H \in C^\infty (M) $ and every pair $m,m' \in M$, the stationary points of the restriction of $ {\mathcal L}^H$ 
to the set of initially cotangent paths for $\pi$ connecting $m$ to $m'$ are cotangent paths for $\pi $.
\item  If for every function $H \in C^\infty (M) $ and every pair $m,m' \in M$, the stationary points of the restriction of $ {\mathcal L}^H$ to the set of 
initially cotangent paths for $\pi$ connecting $m$ to $m'$ are cotangent paths for $\pi$, then the bivector field $\pi$ is weakly foliated.
\item  The bivector field $\pi$ is Poisson if and only if, 
for every function $H \in C^\infty (M) $ and every pair $m,m' \in M$, the stationary points of the restriction of $ {\mathcal L}^H$ to the set
${\tilde P}_{m,m'}(T^*M)$ of paths connecting $m$ to $m'$ are quasi-cotangent paths for $( \pi,H)$.
\end{enumerate}
\end{theorem}

Notice the following straightforward corollary: 

 \begin{corollary}\label{coro:cheminscotangents_si_poisson} 
Let $\pi$ be a bivector field on an open $M$ of $\R^n$. If $ \pi$ is a Poisson structure, then,
 for every functions $H \in \FF{M}$ and every pair of points $m,m'$ in $ M$, the stationary points of the restriction of $ {\mathcal L}^H$ to
 the set $ {\hat P}_{m,m'}(T^*M)$ of initially cotangent paths for $\pi$ connecting $m$ to $m' $ are quasi-cotangent paths for $(\pi,H)$.
 \end{corollary}
This proof of theorem \ref{theo:caract_chemins_cotangents} will be done in several stages.

\begin{proposition}\label{prop:pts_statio}
Let $\pi$ be a bivector field on a manifold $M$.
Chose $\nabla$ a connection on $M$.
\begin{enumerate}
\item The stationary points of the restriction of $ {\mathcal L}^H$ to the set of  paths ${\tilde P}_{m,m'}(T^*M)$ are the paths $\alpha$ with base path $x$ that satisfy $x(0)=m,x(1)=m'$ and the following two equations:
\begin{equation} \label{eq:points_stationnaires}
 \left\{  \begin{array}{rcl}   \nabla_{\frac{\diff x(t)}{\diff t}} \alpha &= & -(K_{x(t)}^H)^{*}  (\alpha(t)) \\ 
           \frac{\diff x(t)}{\diff t} & =&   \left. {\mathcal X}_H \right|_{x(t)}.
           \end{array}    
 \right. 
\end{equation}
with $K^H$ as in lemma \ref{lem:quasi_hami}.
\item 
The stationary points of the restriction of $ {\mathcal L}^H$ to the set of initially cotangent paths for $\pi$ from $m$ to $m'$ ${\hat P}_{m,m'}(T^*M)$ are the paths 
that satisfy (\ref{eq:points_stationnaires}) and $x(0)=m,x(1)=m',\pi_m^\#(a(0))= \left. \frac{\diff x(t)}{\diff t}\right|_{t=0}$.
\end{enumerate}
\end{proposition}
\begin{proof}
A path of the form $(0,\delta_0)$ is automatically tangent to both ${\tilde P}_{m,m'}(T^*M)$ and ${\hat P}_{m,m'}(T^*M)$ if $\delta_0 (0)=0$.
For any stationary point of the restriction on $ {\mathcal L}^H$ of the whole ${\tilde P}_{m,m'}(T^*M)$,
By proposition \ref{prop:calcul_diff}, this implies that for all such $\delta_0 $:
$$ \int_0^1 \langle \delta_0(t),   \left. {\mathcal X}_H \right|_{x(t)} - \frac{\diff x(t)}{\diff t} \rangle dt =0.$$
 This condition is satisfied if and only if the second relation in (\ref{eq:points_stationnaires}) is satisfied.
 
 A path of the form  $(\gamma_0,0)$ is automatically tangent to both ${\tilde P}_{m,m'}(T^*M)$ and ${\hat P}_{m,m'}(T^*M)$ if 
 $\gamma_0$ is zero at  $t=0$ and $t=1$ and $\gamma_0$ has a derivative equal to zero at $t=0$.
By proposition \ref{prop:calcul_diff}, for all such $\delta_0 $:
$$\int_0^1 \left\langle \alpha ,\nabla_{\gamma_0} \X_H + \Tor^\nabla \left(\frac{\diff x (t)}{\diff t} ,\gamma_0   \right) \right\rangle dt
-\int_0^1 \left\langle  \nabla_{\frac{\diff x (t)}{\diff t}} \alpha , \gamma_0   \right\rangle =0 dt.
$$
Since the second relation in (\ref{eq:points_stationnaires}) is satisfied, this amounts to:
\begin{eqnarray}0 &=& \int_0^1 \left\langle \alpha ,\nabla_{\gamma_0} \X_H + \Tor^\nabla\left(\X_H,\gamma_0   \right)\right\rangle dt \\
& & -\int_0^1 \left\langle  \nabla_{\frac{\diff x (t)}{\diff t}} \alpha , \gamma_0   \right\rangle  dt \\
& &= \int_0^1 \left\langle \alpha ,K^H (\gamma_0)\right\rangle 
-\int_0^1 \left\langle  \nabla_{\frac{\diff x (t)}{\diff t}} \alpha , \gamma_0   \right\rangle dt.
\end{eqnarray}
 This condition is satisfied if and only if the first relation in (\ref{eq:points_stationnaires}) is satisfied.
\end{proof}

\begin{corollary}\label{cor:existenceofStationnaryPontsStartingFrom}
Let $\pi$ be a bivector field on a manifold $M$.
If the Hamiltonian field $\X_H$ of $H$ is complete, then:
\begin{enumerate}
\item
For any point $m$ in $M$ and all $\alpha_0 \in T^*_m M $, 
there exists $m'\in M$ such that the restriction of $ {\mathcal L}^H$ to ${\tilde P}_{m,m'}(T^*M)$
admits a stationary point $\alpha (t)  $ with $\alpha(0)=\alpha_0$.
\item 
For any point $m$ in $M$ and all $\alpha_0 \in T^*_m M $
such that $\pi^\#_m(\alpha_0 - \diff_m H)=0$, 
there exists $m'\in M$ such that the restriction of $ {\mathcal L}^H$ to ${\tilde P}_{m,m'}(T^*M)$
admits a stationary point $\alpha$ with $\alpha(0)=\alpha_0$.
\end{enumerate}
\end{corollary}
\begin{proof}
The Hamiltonian field $\X_H$ being complete, the system of differential equations (\ref{eq:points_stationnaires}) admits a solution
for any initial value. For the second of these relations, this is by definition of completeness, while for the first one, it is by  linearity of the remaining equation when $x(t)$ is given.
These solutions are stationary points of the restriction of $ {\mathcal L}^H$ to ${\tilde P}_{m,m'}(T^*M)$
with $\alpha(0)=\alpha_0$ by item \emph{1} in  proposition \ref{prop:pts_statio}.
This proves item \emph{1} of the corollary. If in addition $\pi^\#_m(\alpha_0 - \diff_m H)=0$,
which amounts to $\pi^\#_m(\alpha_0) =\X_H\mid _{m} =\frac{\diff x}{\diff t }\mid _{t=0}$,
the solution of (\ref{eq:points_stationnaires}) starting from $\alpha_0$ is 
a stationary point of the restriction of $ {\mathcal L}^H$ to ${\tilde P}_{m,m'}(T^*M)$
by item \emph{2} in  proposition \ref{prop:pts_statio}.
This proves item \emph{2} of the corollary.
\end{proof}

\begin{lemma}\label{lem:antisym_implique}
Let $M$ be a manifold and $\pi$ a bivector field on $M$.
 For every function $H$ and every torsion free connection $\nabla$
 \begin{equation} \label{eq:equa_diff_dH} 
 \nabla_{\frac{\diff x}{\diff t}} \diff_{x(t)}H  = - {K^H_{x(t)}}^{*} (\diff_{x(t)} H ),
 \end{equation}
with $x(t)$ the integral curve through $m$ of $\X_H$ and
$K^H$ is defined as in proposition  \ref{prop:pts_statio}.
\end{lemma}
\begin{proof}
Let $u$ be an arbitrary vector field. By definition of the torsion: 
  \begin{equation}\label{eq:antisymimplique1} \langle \diff_{x(t)} H , \Tor^\nabla(u,\X_H)\rangle 
= \langle \diff_{x(t)} H , \nabla_{\X_H} u - \nabla_u \X_H -[\X_H,u] \rangle.
\end{equation}
By the definition of brackets of vector fields, we have
$\langle \diff_{x(t)} H , [\X_H,u] \rangle = \X_H(u(H))- u (\X_H H)= \X_H(u(H))$, the last identity
follows from the fact that the $\X_H(H)=0$ by skew-symmetry of $\pi$.
By definition of a connection and its dual, 
$$\X_H(u(H)) =  \langle \nabla_{X_H} u , \diff_{x(t)} H \rangle +  \langle  u , \nabla_{X_H} \diff_{x(t)} H \rangle.$$
So that equation (\ref{eq:antisymimplique1}) amounts to:
$$ \langle \diff_{x(t)} H , \nabla_{u} \X_H \rangle +  \langle  u , \nabla_{X_H} \diff_{x(t)} H \rangle = \langle \diff_{x(t)} H , \Tor^\nabla(u,\X_H)\rangle $$ 
i.e.
$$ \langle (K^H)^{*}(\diff_{x(t)} H) - \nabla_{\X_H} \diff_{x(t)} H  , u \rangle  =0.$$
Since $u$ is arbitrary and $\X_H =  \frac{\diff x}{\diff t}$, the result follows.
\end{proof}

\begin{proposition}\label{prop:jacobi_implique}
 Let $M$ be a manifold and $\pi$ a bivector field on $M$.
 For every function $H$, the Lie derivative of $\pi$ with respect
 to $\X_H$ is given, for every torsion free connection $\nabla$, by:
  \begin{equation}\label{eq:equa_diff_pi} (\nabla_{\frac{\diff x(t)}{\diff t}} \pi)^\#  = K^H_m \circ \pi^\#_m + \pi^\# \circ {K^H_{m}}^{*} -(({\mathcal L}_{\X_H} \pi)_m)^\#, 
   \end{equation}
  with $x(t)$ being the integral curve through $m$ of $\X_H$ 
 and $K^H$ as in lemma \ref{lem:quasi_hami}.
\end{proposition}
We start with a lemma:
\begin{lemma}\label{lem:derivation}
Let $M$ be a manifold. For every vector field $u$,
every multivector field $P$ and every connection $\nabla$:
\begin{equation}\label{eq:deriveeDeLieConnection}
 {\mathcal L}_{u} P = \underline{N} (P)  - \nabla_u P 
\end{equation}
with $N(v):=\nabla_v u+Tor^\nabla(u,v) $  for all vector field $v$ and where $\underline{N} : \wedge^{\bullet  } TM \to  \wedge^{\bullet  } TM $ is the natural extension of $N$ by derivation:
$$ N(v_1 \wedge \dots \wedge  v_k) = N(v_1) \wedge v_2 \wedge \dots \wedge v_k + v_1 \wedge N(v_2) \wedge \dots...$$
\end{lemma}
\begin{proof}
Equation  (\ref{eq:deriveeDeLieConnection}) is true when $P$ is a vector field
(in which case we just recover the definition of the torsion),
$$ {\mathcal L}_{u} P - \underline{N} (P)  + \nabla_u P = [u,P]+\nabla_P u - \nabla_u P =Tor(u,P)$$
and the set of multivector field on which it is true is an algebra (because $ {\mathcal L}_{u}$, $\underline{N}$, $ \nabla_u$
are derivations), so it is valid for all multivector $P$.
\end{proof}
%on %\ref{prop:cheminscotangents_si_poisson}.
We can now prove proposition \ref{prop:jacobi_implique}:

\begin{proof}
 Lemma \ref{lem:derivation}, applied to the vector fields $\X_H$, the bivector field $\pi$ gives:
 $${\mathcal L}_{\X_H} \pi = K^H_* (\pi)  - \nabla_{\X_H}\pi$$
 The formula  $(\underline{N} \pi)^\# = N \pi^\# + \pi^\#{N}^{*}$ being valid for every bivector field and every $1$-$1$-tensor $N$.
 Proposition \ref{prop:jacobi_implique} follows. 
\end{proof}

\begin{corollary}
Let $\pi$ be a bivector field on a manifold $M$. For all $m,m' \in M$ and any stationary point $\alpha$ (with base path $x$) of the restriction on $\tilde{P}_{m,m'}$ or $\hat{P}_{m,m'} $  of $ {\mathcal L}_H $,
the following differential equation
\begin{equation}\label{eq:clef} \nabla_{\frac{\diff x(t)}{\diff t}} c(t) =  K^H_{x(t)}  (c (t)) +  ({\mathcal L}_{\X_H} \pi)^\# (\alpha(t)- \diff_{x(t)} H )
\end{equation}
is satisfied by $ c (t) =  \pi^\#_{x(t)} ( \alpha(t)) - \frac{\diff x(t) }{\diff t} $.
(Here $K^H$ is as in lemma \ref{lem:quasi_hami})
\end{corollary}

Note that $  ({\mathcal L}_{\X_H}\pi)^\# (\diff_{x(t)} H)= 0 $. However, it is desirable to keep this term.

\begin{proof}
It follows from (\ref{eq:equa_diff_pi}) and (\ref{eq:points_stationnaires}) that: 
\begin{eqnarray*}
\nabla_{\frac{\diff x(t)}{\diff t}}  \pi^\#_{x(t)} \alpha(t)
&=& K^H_{x(t)} \circ \pi^\#_{x(t)} (\alpha(t)) +  \pi^\#_{x(t)} \circ {K^H_{x(t)}}^{\perp} (\alpha(t)) \\
&+& ({\mathcal L}_{\X_H}\pi)^\# (\alpha(t))  -  \pi^\#_{x(t)} {K^H_{x(t)}}^{\perp} (\alpha(t))\\
&=& K^H_{x(t)} \circ \pi^\#_{x(t)} (\alpha(t))+ ({\mathcal L}_{\X_H}\pi)^\# (\alpha(t)).
\end{eqnarray*}
It also follows from (\ref{eq:points_stationnaires}) that:
 $$ \nabla_{\frac{\diff x(t)}{\diff t}} =  \pi^\#_{x(t)} \left(\diff_{x(t)} H \right).$$
Hence, according to (\ref{eq:equa_diff_dH}):
 $$\frac{\diff x(t)}{\diff t} \pi^\#_{x(t) } \diff_{x(t)}H  
= K^{H}_{x(t)} \circ \pi^\#_{x(t)} (\diff_{x(t)} H(t))+  ({\mathcal L}_{\X_H}\pi)^\# (\diff_{x(t)} H).
$$
Therefore:
$$ \nabla_{\frac{\diff x(t)}{\diff t}}\frac{\diff x(t)}{\diff t} =  K^{H}_{x(t)} \left(\frac{\diff x(t) }{\diff t}\right)+  ({\mathcal L}_{\X_H}\pi)^\# (\diff_{x(t)} H).$$
By taking the difference of these two equations, one obtains equation~(\ref{eq:clef}).
\end{proof}

We can now prove theorem \ref{theo:caract_chemins_cotangents}.

\begin{proof}
Let us show item \emph{1}.
For every function $H \in C^\infty (M) $ and all $m,m' \in M$, stationary points of the restriction of $ {\mathcal L}^H$ 
to the subset ${\hat P}_{m,m'}(T^*M)$ of initially cotangent paths for $\pi$ connecting $m$ to $m'$, are such that:
 $$c(t) = \pi^\#_{x(t)} \left((\alpha(t)) - \frac{\diff x(t) }{\diff t}\right),$$
 satisfies (\ref{eq:clef}).
Now by proposition \ref{prop:11tensor}, if the bivector field $\pi$ is foliated, then
there exists a smooth map $t \to (C_H)_t$ with $(C_H)_t \in End(T_{x(t)}M)$ such that
$({\mathcal L}_{\X_H} \pi)^\# = C_H \circ \pi^\# $ and (\ref{eq:clef})
becomes:
 $$ \frac{\diff (c(t))}{\diff t} =  K^{H}_{x(t)}  (c (t)) +  (C_H)_{t}  \pi^\#_{x(t)} \left(  (\alpha(t)) - {\diff_{x(t)} H} \right) 
 =   M_{t}(c(t)) $$
 with $M_{t}:=  K^{H}_{x(t)}  + (C_H)_t$.
 This equation is linear, and by definition of initially cotangent paths, $c(0)=0$ so that $c(t)=0$ for all $t$,
 which means that the stationary point $\alpha$ is a cotangent path for $\pi$.

Let us show item \emph{2}.
Let $\alpha$ be a stationary point of $ {\mathcal L}^H$ restricted to ${\hat P}_{m,m'}(T^*M)$.
By assumption it is a cotangent path for $\pi$, 
so that the path $c$ defined as in (\ref{eq:clef}) is equal to $0$. Equation 
(\ref{eq:clef}) then implies that, for every stationary point $\alpha$  of $ {\mathcal L}^H$ restricted to ${\hat P}_{m,m'}(T^*M)$, the following equation holds:
 $$ ({\mathcal L}_{\X_H} \pi)^\# (\alpha(t)- \diff_{x(t)} H ) =0 ,$$
 where $x$ is, as usual, the base path.
 
Now for all  $m\in M $, and all function $H $  defined in a neighborhood of $m$,
there exists a function with compact support, that we still denote by $H$, that coincides with $H$
in a  neighborhood of $m$.
The Hamiltonian vector field $\X_H$ is complete, 
so that for all $\alpha_0 \in T^*_m M$ with $\pi^\#_m(\alpha_0) = (\X_H){_{\vert m}} $, i.e. such that:
 \begin{equation}\label{eq:pidies} 
 \pi^\#_m(\alpha_0- \diff_{m} H ) =0,
\end{equation}
 there exists by corollary \ref{cor:existenceofStationnaryPontsStartingFrom} a stationary point  
 $\alpha$ of the restriction of $ {\mathcal L}^H$ to ${\hat P}_{m,m'}(T^*M)$,
with $m'=\phi_1(m)$ starting from $\alpha_0$. 
 As we just saw, this stationary point $\alpha$ satisfies
 \begin{equation}\label{eq:pidiestwo} 
({\mathcal L}_{\X_H} \pi )^\# ( \alpha- \diff_{x(t)} H ) =0.
\end{equation}
 Applied to $t=0$, this identity amounts to $({\mathcal L}_{\X_H} \pi)^\# ( \alpha_0- \diff_{m} H )$.
 Hence for every $\alpha_0$ such that equation \ref{eq:pidies} holds, equation \ref{eq:pidiestwo} also holds, which means that the kernel of  $\pi^\#_m$ is contained in the kernel of $(({\mathcal L}_{\X_H} \pi)_m)^\# $,
 or, by skew-symmetry, that the image of $(({\mathcal L}_{\X_H} \pi)_m)^\# $ is contained in the image of $ \pi^\#$.
 Since this is valid for all $m \in M$ and all function $H$, proposition \ref{foliated=imageweakly} implies that $\pi$ is weakly foliated.

Let us show item \emph{3}. 
Let $H$ be a function in $C^\infty (M) $ and let $\alpha $ be a stationary point with base path of the restriction of $ {\mathcal L}^H$
to  the set of paths relating $m$ and $m'$.
Since $\pi$ is Poisson, the Lie derivative of $\pi$ in the direction of $\X_H$ vanishes, and Equation (\ref{eq:clef}) becomes
 $$  \frac{\diff (c(t))}{\diff t} =  K^H_{x(t)}  (c (t)) \textrm{~~~~with~~~~} c(t)=\pi^\#(\alpha(t))-\frac{\diff x(t) }{\diff t}.$$
By lemma \ref{lem:quasi_hami}, paths that satisfy this equation are quasi-cotangent paths for $(\pi,H)$.

Conversely,
for all  $m\in M $ and $\alpha_0 \in T^*_mM$, and all function $H $  defined in a neighborhood of $m$,
there exists a function with compact support, that we still denote by $H$, that coincides with $H$
in a  neighborhood of $m$.
The Hamiltonian vector field $\X_H$ is complete,   
so that  there exists by Corollary \ref{cor:existenceofStationnaryPontsStartingFrom} a stationary point  
$\alpha$ of the restriction of $ {\mathcal L}^H$ to ${\tilde P}_{m,m'}(T^*M)$,
with $m'=\phi_1(m)$ starting from $\alpha_0$. 
By assumption, $\alpha$ is quasi-cotangent for $(\pi,H)$. Lemma \ref{lem:quasi_hami} implies that it satisfies:
 $$ \frac{\diff (c(t))}{\diff t} =  K^{H}_{x(t)}  (c (t)) \textrm{~~~~with~~~~} c(t)=\alpha(t)-\diff_{x(t)}H,$$ 
 this is compatible with equation (\ref{eq:clef}) if and only if:
  $$ ({\mathcal L}_{\X_H} \pi)^\# (\alpha(t)- \diff_{x(t)} H )  =0.$$
  Applied to $t=0$, this identity amounts to $({\mathcal L}_{\X_H} \pi)^\# ( \alpha_0- \diff_{m} H )$.
  Since $\alpha_0$ is arbitrary, this implies that ${\mathcal L}_{\X_H}\pi=0$. Since $H$ is also arbitrary (at least in a neighborhood of $m$),
  $\pi$ needs to be Poisson in a neighborhood of $m$. Since $m$ is arbitrary, the bivector field $\pi$ is Poisson.
\end{proof}

\subsection{Several counterexamples}
\label{sec:counterexamples}

We explain why  theorem \ref{theo:caract_chemins_cotangents} is the best outcome we can hope for,
by giving counterexamples to its attempted reciprocal.

\noindent {\bf I. The implication stated in item \emph{1)} of theorem \ref{theo:caract_chemins_cotangents} can not be inverted.} 

Consider the bivector field (\ref{eq:explicitNotFoliated}) on $M={\mathbb R}^4$ of example \ref{ex:WeaklyNotStrongly} for $i=1$
which is not foliated (although it is weakly foliated).
We claim that for every function $H \in C^\infty (M) $ and every pair $m,m' \in M$, the stationary points of the restriction of $ {\mathcal L}^H$ 
to the set of initially cotangent paths for $\pi$ connecting $m$ to $m'$ are cotangent paths for $\pi$.

Specializing relations (\ref{eq:points_stationnaires}) to the canonical connection on $M={\mathbb R}^4$, one sees that stationary paths for a given function $H$ 
of the restriction of $ {\mathcal L}^H$ to the set of initially cotangent paths for $\pi$ connecting $m$ to $m'$
are the solutions $\alpha(t)=(x(t),a(t)) \in {\mathbb R}^8 \simeq  T^* {\mathbb R}^4$ of the following differential equations:
 $$
 \left\{  \begin{array}{rcl}  \frac{\diff a(t)}{\diff t} &= & -{K^H_{x(t)}}^{*}  (a(t)) \\ 
           \frac{\diff x(t)}{\diff t} & =&   \left. {\mathcal X}_H \right|_{x(t)}.
           \end{array}    
 \right. 
$$
If $m$ is a point of the form $m=(0,0,u,v)$, then $\pi_m$ is zero at this point. Therefore
 either a stationary path of the restriction of $ {\mathcal L}^H$ is for all $t$ at points of the form $m=(0,0,u,v)$, or it never intersects such a point.
If $m$ does not belong to this set, $\pi$ is then foliated, and therefore the stationary paths of the restriction of $ {\mathcal L}^H$ 
to the set of initially cotangent paths for $\pi$ connecting $m$ to $m'$ are cotangent paths for $\pi$ by item 1) in Theorem
\ref{theo:caract_chemins_cotangents}. We therefore an example for which all the stationary points are cotangent path although the bivector field is not foliated.

\noindent {\bf II. The implication stated in item \emph{2)} of theorem \ref{theo:caract_chemins_cotangents} can not be inverted.} 

Consider the bivector field (\ref{eq:explicitNotFoliated}) on $M={\mathbb R}^4$ of example \ref{ex:WeaklyNotStrongly} for $i=0$. Again, this bivector field is not foliated, although it is weakly foliated.
Choose $H(x,u,y,v)=u$.
The path $\alpha(t)=(a(t),x(t))$ in $T^*M \simeq M \times M \simeq {\mathbb R}^4 \times {\mathbb R}^4$ 
given by
$$a(t)=(0,1,1,0) \hbox{ and } x(t)=(t,0,0,0)$$
has a base path $x(t)$ admitting $m=(0,0,0,0) \in M$ as starting point and  $m'=(1,0,0,0) \in M$
as ending point.

This path is initially cotangent for $\pi$, since 
$$\pi_{x(0)}^\#(0,1,1,0) =(1,0,0,0)= \left. \frac{\diff x(t)}{\diff t}\right|_{t=0} .$$
It is however not a a cotangent path for $\pi$, because this identity is not valid any more for $t \neq 0$. 

Now, we check that $\alpha$ is a stationary point for the restriction of   ${\mathcal L}^H $ to ${\tilde P}_{m,m'}(T^*M)$. 
 By proposition \ref{prop:pts_statio}, it suffices to verify that equation (\ref{eq:points_stationnaires}) are satisfied,
 which is an obvious verification in this case. First $x(t)=\X_H =\frac{\partial }{\partial x} $. Now, since $\X_H = \frac{\partial }{\partial x}$, 
 we have $K^H=0$ when computed with respect to the canonical connection on ${\mathbb R}^4$, while $ \frac{\diff a(t)}{\diff t}= 0$ so that the second relation in (\ref{eq:points_stationnaires}) 
 holds as well. Hence we have an example of a weakly foliated bivector field for which there is a function $H$ and a points $m,m'\in M$ such that the restriction of $ {\mathcal L}^H$ to the set
${\tilde P}_{m,m'}(T^*M)$ admits a stationary points which is not a cotangent path.
\section{Relation with Poisson $\sigma$-models}
 \label{sec:relation}
 
 In \cite{strobelklim}, Klim{\v{c}}{\'{\i}}k and Strobl introduce a functional whose stationary points correspond to Lie algebroid morphisms.
We would like to relate our functional to theirs. A technical difficulty comes from the fact that they
work with loops while we work with paths.
Adapted to paths, the functional that appear in (1) of \cite{KosmannSchwarzbachYvette} is defined as follows.
Let $M$ be a manifold equipped with a bivector field $\pi$, a priori not Poisson.
Consider $\beta  $ a vector bundle morphism from $T(I^2) $ to $T^*M $
with $I = [0,1]$. Concretely, $\beta$ is defined by a map $X: I^2 \to M$
and a pair of maps from $I^2$ to $T^*M $ above $X$, namely $\beta \left(\pp{}{ t}\right) $
and $  \beta \left(\pp{}{y}\right) $, with $t$ and $y$ being the canonical coordinates of the square $I^2$.
The functional of \cite{strobelklim} is then:

  \begin{eqnarray*}\label{eq:KSLagrangien} {\mathcal L}^{KS} &=&  
 \int_{I^2} \left(  \left\langle \beta \left(\pp{}{ t}\right), \pp{X}{y}  \right\rangle -  \left\langle \beta \left(\pp{}{ y}\right), \pp{X}{t}   \right\rangle \right) \diff t \diff y\\
  &+& \int_{I^2} \left\langle  \beta \left(\pp{}{ t}\right) \wedge \beta \left(\pp{}{ y}\right),\pi_{X(y,t)} \right\rangle   \diff t \diff y.
 \end{eqnarray*}
Out of every map  $\alpha : I \to T^* M $ with base path $x $ we can, given a function $H$ whose Hamiltonian flow 
is well-defined at time $1$ (at least at all points of the curve $x$) associate a vector bundle morphism $\tilde{\alpha}$ from 
 $ T(I^2) $ to $T^* M$ by  $ X(t,y) = \phi_y (x(t))$ where $\phi$ is the flow of $\X_H$
 and the vector bundle morphism induced by:
  $$ \tilde{\alpha}\left(\pp{}{ y}\right)  = T\phi_{-y}^* \left(\alpha(t)\right) \hbox{ ~~and~~ }  \tilde{\alpha} \left(\pp{}{ t}\right) = \diff_{X(t,y)} H .$$
  (By construction, both $\tilde{\alpha}( \pp{}{ t}) $ and $\tilde{\alpha}(\pp{}{ y}) $,
  defined as previous, belong to the cotangent space $T_{X(t,y) }^* M $ for all $(y,t) \in I^2$,
  which justifies the definition of $\tilde{\alpha}$.)
 
\begin{proposition}
Let $M$ be a manifold equipped with a bivector field $\pi$.
For all function $H$ whose Hamiltonian flow is defined at time $1$
and all $\alpha : I \to T^*M$, the following relation holds: 
 $$ {\mathcal L}^H (\alpha) = {\mathcal L}^{KS}(\tilde{\alpha}) ,$$
 with $\tilde{\alpha} $ being defined as above.
\end{proposition}
\begin{proof}
According to (\ref{eq:KSLagrangien}),  the quantity ${\mathcal L}^{KS}(\tilde{\alpha}) $ is the alternate sum of three integrals:
$$ \begin{array}{l}\int_{I^2} \langle \tilde{\alpha} \left(\pp{}{ t}\right), \pp{X}{y}  \rangle\diff t \diff y,\int_{I^2} \langle \tilde{\alpha} \left(\pp{}{ y}\right),
 \pp{X}{t}  \rangle\diff t \diff y, \\ \int_{I^2}  \langle \tilde{\alpha} \left(\pp{}{ t}\right) \wedge \tilde{\alpha} \left(\pp{}{ y}\right),\pi_{X(t,y)}   \rangle\diff t \diff y  \end{array}$$
We claim that the integrand of the first term vanishes, while the integrands of the second and third terms does not depend on the variable~$y$. 
Let us prove these three points.
First, by construction $ X(t,y) = \phi_y (x(t))$, so that $\pp{X}{y} = \left. \X_H\right|_{\phi_y (x(t))} $.
By definition also, $\tilde{\alpha} (\pp{}{ t}) = \diff_{X(t,y)} H $
so that 
$$
 \left\langle \tilde{\alpha} \left(\pp{}{ t}\right), \pp{X}{y} \right\rangle=\langle \diff_{X(t,y)} H, \left. \X_H\right|_{\phi_y (x(t))}  \rangle = 0, 
 $$
by skew-symmetry of $\pi$. This proves the first point, which immediately implies:
\begin{equation}
 \label{eq:KSLagrangian0} 
  \int_{I^2}\left\langle \tilde{\alpha} \left(\pp{}{ t}\right), \pp{X}{y}  \right\rangle\diff t \diff y =0.
\end{equation}

For the second point, it suffices to see that for all $t\in I$ and all $y \in I$:
 \begin{eqnarray*} \left\langle \tilde{\alpha} \left(\pp{}{ y}\right), \pp{X}{t}  \right\rangle &=&   \left\langle  T\phi_{-y}^* \alpha(t)  , T\phi_y ( \frac{\diff x(t)}{\diff t}) \right\rangle \\
  &=& \left\langle  \alpha(t)  , \frac{\diff x(t)}{\diff t} \right\rangle.
 \end{eqnarray*}
This quantity does not depend on $y$, which, in turn, gives the relation:
\begin{equation}
 \label{eq:KSLagrangian1}
  \int_{I^2} \left\langle \tilde{\alpha} \left(\pp{}{ y}\right), \pp{X}{t}  \right\rangle\diff t \diff y =  \int_{I} \left\langle  \alpha(t)  , \frac{\diff x(t)}{\diff t}  \right\rangle\diff t.
\end{equation}

For the third point, it suffices to see that for all $t\in I$ and all $y \in I$, by definition:
 \begin{eqnarray*} \left\langle  \tilde{\alpha} \left(\pp{}{ t}\right) \wedge \tilde{\alpha} \left(\pp{}{ y}\right),\pi_{X(t,y)}   \right\rangle &=&  
 \left\langle  \diff_{X(t,y)} H \wedge T\phi_{-y}^* \alpha(t) ,   \pi_{X(t,y)} \right\rangle \\
  &=& \left\langle T\phi_{-y}^* \alpha(t),  \left. \X_H\right|_{\phi_y (x(t))} \right\rangle.
 \end{eqnarray*}
A vector field being always preserved by its flow, we have:
$$ \left. \X_H\right|_{\phi_y (x(t))} = T_{\phi_0(x(t))} \phi_y \left( \left. \X_H\right|_{\phi_0(x(t))} \right) =T_{x(t)} \phi_y \left( \left. \X_H\right|_{x(t)} \right) ,$$
which implies:
 \begin{eqnarray*}\left\langle  \tilde{\alpha} \left(\pp{}{ t}\right) \wedge \tilde{\alpha} \left(\pp{}{ y}\right),\pi_{X(y,t)}   \right\rangle
&=& \left\langle T\phi_{-y}^* \alpha(t), T_{x(t)} \phi_y \left( \left. \X_H\right|_{x(t)} \right)  \right\rangle \\ & =
& \left\langle \alpha(t), \left( \left. \X_H\right|_{x(t)} \right)  \right\rangle.  
\end{eqnarray*}
This quantity does not depend on $y$, which in turn gives the relation:
\begin{equation}
 \label{eq:KSLagrangian2}
  \int_{I^2}  \left\langle \tilde{\alpha} \left(\pp{}{ t}\right) \wedge \tilde{\alpha} \left(\pp{}{ y}\right),\pi_{X(t,y)}   \right\rangle\diff t \diff y 
  = \int_I  \left\langle \alpha(t), \left( \left. \X_H\right|_{x(t)} \right)  \right\rangle \diff t.
\end{equation}
The alternate sum of the right hand sides of (\ref{eq:KSLagrangian0}-\ref{eq:KSLagrangian1}-\ref{eq:KSLagrangian2}) gives ${\mathcal L}^H (\alpha)$,
while the alternate sum of the left hand sides of these three equations gives ${ \mathcal L}^{KS}(\tilde{\alpha}) $. This completes the proof.
 \end{proof}

  % the introduction
%\newpage

 %references
 
\end{document}